\documentclass[final,12pt,a4paper]{article}

\usepackage{showkeys}
\usepackage[english,francais]{babel} 
\usepackage{ucs}
\usepackage[utf8]{inputenc}

\RequirePackage{color} 

 \usepackage{hyperref}
 \usepackage{graphicx}
 \DeclareGraphicsExtensions{.eps, .ps, .pstex} 
\RequirePackage{amsmath,amssymb,amsfonts,bbm,latexsym,mathrsfs}

\newcommand{\ome}{\ov{\mathbbm{e}}}
\newcommand{\oz}{\ov{Z}}

 \newcommand{\X}{{\rm X}}
\newcommand{\s}{{\rm S}}
\newcommand{\pmc}{{\rm P}\mathbb{M}}

\newcommand{\dgr}{d^{\rm gr}} \newcommand{\dgh}{\mathtt{d}_{{\rm GH}}}

 \newcommand{\bq}{{\bf q}}

\newcommand{\bl}{{\bf l}} 
 \newcommand{\bt}{{\bf t}}

 \newcommand{\bT}{{\bf T}}
 
\newcommand{\bQ}{{\bf Q}}

  \newcommand{\R}{\mathbb{R}}
\newcommand{\N}{\mathbb{N}} \newcommand{\M}{\mathbb{M}}

\newcommand{\TT} { {\cal T }}

\def\build#1_#2^#3{\mathrel{ \mathop{\kern 0pt#1}\limits_{#2}^{#3}}}

\def\cq{$\hfill \square$}

\def\e{{\cal E}}
\def\s{{\cal S}}

\def\v{{\cal V}}

\def\eps{\varepsilon}

\def\ov{\overline}

\def\rems{\noindent{\bf Remarks. }}

\def\proof{\noindent{\bf Proof. }}

\def\diam{{\rm diam\,}}

\newtheorem{thm}{Theorem}
\newtheorem{lmm}{Lemma}
\newtheorem{prp}{Proposition}
\newtheorem{defn}{Definition}

\usepackage{fancyhdr}

\usepackage[outerbars]{changebar} 

\def\miermont{{\href{http://mahery.math.u-psud.fr/~miermont/}{Gr\'egory
      Miermont}}}
 
\def\LPMA{{\href{http://www.proba.jussieu.fr}{Laboratoire de
      Probabilit\'es et Mod\`eles Al\'eatoires}}}
\def\cnrs{{\href{http://www.cnrs.fr/}{CNRS}}}
\def\FSMP{{\href{http://www.sciencesmaths-paris.fr}{Fondation des
      Sciences Math\'ematiques de Paris}}}

\title{On the sphericity of scaling limits of random planar
  quadrangulations}

\author{ \miermont \thanks{\LPMA, 
Universit\'e Pierre et Marie Curie, 175 rue du
    Chevaleret, F-75013 Paris. {\tt
      Gregory.Miermont@math.u-psud.fr}. This research is supported by
    the \cnrs\ and the \FSMP.}  }

\begin{document}

\selectlanguage{english}

\maketitle

%\tableofcontents

\begin{abstract}
We give a new proof of a theorem by Le Gall \& Paulin, showing that
scaling limits of random planar quadrangulations are homeomorphic to
the $2$-sphere. The main geometric tool is a reinforcement of the
notion of Gromov-Hausdorff convergence, called $1$-regular
convergence, that preserves topological properties of metric surfaces.
\end{abstract}

\section{Introduction}\label{sec:introduction}

A planar map is a combinatorial embedding of a connected graph into
the 2-dimensional sphere. Random planar map have drawn much attention
in the recent probability literature due to mathematical physics
motivations \cite{ADJ} and a powerful encoding of planar maps in terms
of labeled trees due to Schaeffer \cite{schaeffer98,CSise}. In turn,
scaling limits of labeled trees are well-understood thanks to the
works of Aldous, Le Gall and others
\cite{aldouscrt93,legall99,legall05}. Using this line of reasoning,
many results have been obtained on the geometric aspects of large
random quadrangulations (where faces all have degree $4$), and other
families of maps. Le Gall \cite{legall06} showed in particular that
scaling limits of random quadrangulations are homeomorphic to the
Brownian map introduced by Marckert \& Mokkadem \cite{MM05}, and Le
Gall \& Paulin \cite{lgp} showed that the topology of the latter is
that of the 2-dimensional sphere, hence giving a mathematical content
to the claim made by physicists that summing over large random
quadrangulations amount to integrating with respect to a (still
ill-defined) measure over surfaces.

The aim of this note is to give an alternative proof of Le Gall \&
Paulin's result.  We still strongly rely on the results established by
Le Gall \cite{legall06}, but use very different methods from those of
\cite{lgp}, where the reasoning uses geodesic laminations and a
theorem due to Moore on the topology of quotients of the sphere. We
feel that our approach is somewhat more economic, as it only needs
certain estimates from \cite{legall06} and not the technical
statements of \cite[Lemmas 3.1, 3.2]{lgp} that are necessary to apply
Moore's theorem. On the other hand, this is at the cost of checking
that quadrangulations are close to being path metric spaces, which is
quite intuitive but needs justification (see definitions below). Our
main geometric tool is a reinforcement of Hausdorff convergence,
called $1$-regular convergence and introduced by Whyburn, and which
has the property of conserving the topology of surfaces. We will see
that random planar quadrangulations converge $1$-regularly, therefore
entailing that their limits are of the same topological nature. In the
case, considered in this paper, of surfaces with the topology of the
sphere, the $1$-regularity property is equivalent to \cite[Corollary
  1]{lgp}, stating that there are no small loops separating large
random quadrangulations into two large parts. We prove this by a
direct argument rather than obtaining it as a consequence of the
theorem.

The basic notations are the following. We let $\bQ_n$ be the set of
rooted\footnote{Which means that one oriented edge of the
  quadrangulation is distinguished as the root} quadrangulations of
the sphere with $n$ faces, which is a finite set of cardinality
$2\cdot 3^n(2n)!/(n!(n+2)!)$, see \cite{CSise}. We let $\bq_n$ be a random
variable picked uniformly in $\bQ_n$, and endow the set $V(\bq_n)$ of
its vertices with the usual graph distance $\dgr_n$,
i.e.\ $\dgr_n(x,y)$ is the length of a minimal (geodesic) chain of
edges going from $x$ to $y$.

We briefly give the crucial definitions on the Gromov-Hausdorff
topology, referring the interested reader to \cite{burago01} for more
details. The isometry class $[X,d]$ of the metric space $(X,d)$ is the
collection of all metric spaces isometric to $(X,d)$.  We let $\M$ be
the set of isometry-equivalence classes of {\em compact} metric
spaces. The latter is endowed with the Gromov-Hausdorff distance
$\dgh$, where $\dgh(\X,\X')$ is defined as the least $r>0$ such that
there exist a metric space $(Z,\delta)$ and subsets $X,X'\subset Z$
such that $[X,\delta]=\X,[X',\delta]=\X'$, and such that the Hausdorff
distance between $X$ and $X'$ in $(Z,\delta)$ is less than or equal to
$r$. This turns $\M$ into a complete separable metric space, see
\cite{evpiwin} (this article focuses on compact $\R$-trees, which form
a closed subspace of $\M$, but the proofs apply without change to
$\M$).

\begin{thm}[\cite{lgp}]\label{T1}
A limit in distribution of $[V(\bq_n),n^{-1/4}\dgr_n]$ for the
Gromov-Hausdorff topology, where $n\to\infty$ along some subsequence,
is homeomorphic to the $2$-sphere.
\end{thm}

\rems
\noindent$\bullet$ One of the main open questions in the topic of
scaling limits of random quadrangulations is to uniquely characterize
the limit, i.e.\ to get rid of the somewhat annoying ``along some
subsequence'' in the previous statement.

\medskip

\noindent$\bullet$ To be perfectly accurate, Le Gall \& Paulin showed
the same result for uniform $2k$-angulations (maps with degree-$2k$
faces) with $n$ faces. Our methods also apply in this setting (and
possibly to more general families of maps), but we will restrict
ourselves to the case of quadrangulations for the sake of brevity.

\medskip

\noindent$\bullet$ In the work in preparation \cite{miertopog}, we
provide a generalization of this result to higher genera, in the
framework of Boltzmann-Gibbs distributions on quadrangulations rather
than uniform laws.

\bigskip

As we are quite strongly relying on Le Gall's results in
\cite{legall06}, we will mainly focus on the new aspects of our
approach. As a consequence, this paper contains two statements whose
proofs will not be detailed (Proposition
\ref{sec:estim-lengths-geod-2} and Lemma
\ref{sec:estim-lengths-geod-1}), because they are implicit in
\cite{legall06} and follow directly from the arguments therein, and
also because their accurate proof would need a space-consuming
introduction to continuum tree and snake formalisms. Taking them for
granted, the proofs should in a large part be accessible to readers
with no particular acquaintance with continuum trees or Schaeffer's
bijection.

\section{Gromov-Hausdorff convergence and regularity}\label{sec:grom-hausd-conv}

We say that a metric space $(X,d)$ is a path metric space if every two
points $x,y\in X$ can be joined by a path isometric to a real segment
(with length $d(x,y)$). We let $\pmc$ be the set of isometry classes
of compact path metric spaces, and the latter is a closed subspace of
$(\M,\dgh)$, see \cite[Theorem 7.5.1]{burago01}. One of the main tools
needed in this article is a notion that reinforces the convergence in
the metric space $(\pmc,\dgh)$, which was introduced by Whyburn in
1935 and was extensively studied in the years 1940's. Our main source
is Begle \cite{begle}.

\begin{defn}\label{sec:grom-hausd-conv-2}
Let $(\X_n,n\geq 1)$ be a sequence of spaces in $\pmc$ converging to a
limit $\X$. We say that $\X_n$ converges $1$-regularly to $\X$ if for
every $\eps>0$, one can find $\delta,N>0$ such that for all $n\geq N$,
every loop in $X_n$ with diameter $\leq \delta$ is homotopic to $0$ in
its $\eps$-neighborhood.
\end{defn}

There are a couple of slight differences between this definition and
that in \cite{begle}. In the latter reference, the setting is that
$X_n$ are compact subsets of a common compact space, converging in the
Hausdorff sense to a limiting set $X$. This is not restrictive as
Gromov-Hausdorff convergence entails Hausdorff convergence of
representative spaces in a common compact space, see for instance
\cite[Lemma A.1]{GPW06}. It is also assumed in the definition of
$1$-regular convergence that for every $\eps>0$, there exists
$\delta,N>0$ such that any two point that lie at distance $\leq
\delta$ are in a connected subset of $\X_n$ of diameter $\leq \eps$,
but this condition is tautologically satisfied for path metric
spaces. Last, the definition in \cite{begle} is stated in terms of
homology, so our definition in terms of homotopy is in fact stronger.

The following theorem is due to Whyburn, see \cite[Theorem 6]{begle}
and comments before.

\begin{thm}\label{sec:grom-hausd-conv-1}
Let $(\X_n,n\geq 1)$ be a sequence of elements of $\pmc$ that are all
homeomorphic to $\mathbb{S}^2$. Assume that $\X_n$ converges to $\X$
for the Gromov-Hausdorff distance, where $\X$ is not reduced to a
point, and that the convergence is $1$-regular. Then $\X$ is
homeomorphic to $\mathbb{S}^2$ as well.
\end{thm}

\section{Quadrangulations}\label{sec:quadrangulations}

Rooted quadrangulations are maps whose faces all have degree $4$, and their
set is denoted by $\bQ:=\bigcup_{n\geq 1}\bQ_n$ with the notations of
the Introduction. For $\bq\in \bQ$ we let $V(\bq),E(\bq),F(\bq)$ be
the set of vertices, edges and faces of $\bq$, and denote by
$\dgr_\bq$ the graph distance on $V(\bq)$.

\subsection{A metric surface representation}\label{sec:turn-quadr-into}

One of the issues that must be addressed in order to apply Theorem
\ref{sec:grom-hausd-conv-1} is that the metric space
$[V(\bq),\dgr_\bq]$ is not a surface, rather, it is a finite metric
space. We take care of this by constructing a particular graphical
representative of $\bq$ which is a path metric space whose restriction
to the vertices of the graph is isometric to $(V(\bq),\dgr_\bq)$.

Let $(X_f,d_f),f\in F(\bq)$ be copies of the emptied unit cube ``with
bottom removed''
$$X_f=[0,1]^3\setminus (0,1)^2\times [0,1)\, ,$$ endowed with the
  intrinsic metric $d_f$ inherited from the Euclidean metric
  (i.e.\ the distance between two points of $X_f$ is the minimal
  Euclidean length {\bf of a path in $X_f$}).  Obviously each
  $(X_f,d_f)$ is a path metric space homeomorphic to a closed disk of
  $\R^2$. For each face $f\in F(\bq)$, we label the four incident
  half-edges turning counterclockwise as $(e_1,e_2,e_3,e_4)$, where
  the labeling is arbitrary among the $4$ possible labelings
  preserving the cyclic order. Then define
$$\begin{array}{lll}
c_{e_1}(t)=(t,0,0)_f&\, ,&\qquad 0\leq t\leq
1\\ c_{e_2}(t)=(1,t,0)_f&\, ,&\qquad 0\leq t\leq
1\\ c_{e_3}(t)=(1-t,1,0)_f&\, ,&\qquad 0\leq t\leq
1\\ c_{e_4}(t)=(0,1-t,0)_f&\, ,&\qquad 0\leq t\leq 1\, .
\end{array}$$
In these notations, we keep the subscript $f$ to differentiate points
of different spaces $X_f$. In this way, for every $e\in E(\bq)$, we
have defined a path $c_e$ of length $1$ which goes along one of the
four edges of the boundary $\partial X_f=([0,1]^2\setminus
(0,1)^2)\times \{0\}$, where $f$ is the face incident to $e$.

We then define an equivalence relation $\sim$ on the disjoint union
$\amalg_{f\in F(\bq)}X_f$, as the coarsest equivalence relation such
that for every $e\in E(\bq)$, and every $t\in[0,1]$, we have
$c_e(t)\sim c_{\ov{e}}(1-t)$.  By identifying points of the same
class, we glue the boundaries of the spaces $X_f$ together in a way
that is consistent with the map structure. More precisely, the
topological quotient $\s_\bq:=\amalg_{f\in F(\bq)}X_f/\sim$ is a
$2$-dimensional cell complex whose $1$-skeleton $\e_\bq$ is a graph
representation of $\bq$, and where the faces are the interiors of the
spaces $X_f$. In particular, $\s_\bq$ is homeomorphic to
$\mathbb{S}^2$. We let $\v_\bq$ be the $0$-skeleton of this complex,
i.e.\ the vertices of the graph. We call the $1$-cells and $0$-cells
of $\e_\bq$ and $\v_\bq$ the edges and vertices of $\s_\bq$.

We next endow the disjoint union $\amalg_{f\in F(\bq)} X_f$ with the
largest pseudo-metric $D_\bq$ that is compatible with $d_f,f\in
F(\bq)$ and with $\sim$, in the sense that $D_\bq(x,y)\leq d_f(x,y)$
for $x,y\in X_f$, and $D_\bq(x,y)=0$ for $x\sim y$. Therefore, the
function $D_\bq:\amalg_{f\in F(\bq)} X_f\times \amalg_{f\in F(\bq)}
X_f\to \R_+$ is compatible with the equivalence relation, and its
quotient mapping $d_\bq$ defines a pseudo-metric on the quotient space
$\s_\bq$.

\begin{prp}\label{sec:turn-quadr-into-1}
The space $(\s_\bq,d_\bq)$ is a path metric space homeomorphic to
$\mathbb{S}^2$. Moreover, the restriction of $\s_\bq$ to the set
$\v_\bq$ is isometric to $(V(\bq),\dgr_\bq)$, and any geodesic path in
$\s_\bq$ between two elements of $\v_\bq$ is a concatenation of edges
of $\e_\bq$. Last,
$$\dgh([V(\bq),\dgr_\bq],[\s_\bq,d_\bq])\leq 3\, .$$
\end{prp}

\proof What we first have to check is that $d_\bq$ is a true metric on
$\s_\bq$, i.e.\ it separates points. To see this, we use the fact
\cite[Theorem 3.1.27]{burago01} that $D_\bq$ admits the constructive
expression:
$$D_\bq(a,b)=\inf\left\{ \sum_{i=0}^nd(x_i,y_i):n\geq
0,x_0=a,y_n=b,y_i\sim x_{i+1}\right\}\, ,$$ where we have set
$d(x,y)=d_f(x,y)$ if $x,y\in X_f$ for some $f$, and $d(x,y)=\infty$
otherwise. It follows that for $a\in X_f\setminus \partial X_f$, and
for $b\neq a$, $D_\bq(a,b)>\min(d(a,b),d_f(a,\partial X_f))>0$, so $a$
and $b$ are separated.

It remains to treat the case $a\in \partial X_f$ for some $f$.  The
crucial observation is that a shortest path in $X_f$ between two
points of $\partial X_f$ is entirely contained in $\partial
X_f$. Therefore, the distance $D_\bq(a,b)$ is always larger than the
length of a path with values in the edges $\amalg \partial X_f/\sim$
of $\s_\bq$, where all edges have total length $1$. In particular,
points in distinct classes are at positive distance. One deduces that
$d_\bq$ is a true distance on $\s_\bq$, and by the compactness of the
latter, $(\s_\bq,d_\bq)$ is homeomorphic to $\mathbb{S}^2$
\cite[Exercise 3.1.14]{burago01}.

From this same observation, we obtain that a shortest path between
vertices of $\s_\bq$ is a shortest path of edges, i.e.\ is the
geodesic distance for the (combinatorial) graph distance. Thus,
$(\v_\bq,d_\bq)$ is indeed isometric to $(V(\bq),\dgr_\bq)$. The last
statement follows immediately from this and the fact that
$\diam(X_f,d_f)\leq 3$, entailing that $\v_\bq$ is $3$-dense in
$(\s_\bq,d_\bq)$, i.e.\ its $3$-neighborhood in $(\s_\bq,d_\bq)$
equals $\s_\bq$. \cq

%Note that if $\gamma$ is a loop of diameter $\leq D$ in $S_\bq$, then
%we can find a simple loop $\gamma'$ of diameter $\leq D+4$ in $S_\bq$
%that takes its values only in edges of $S_\bq$, i.e.\ only in the
%points of $\amalg_{f\in F(\bq)}\partial X_f/\sim$. Such a simple path
%is obtained as the boundary of the faces of $S_\bq$ that contain some
%point of $\gamma$. Consequently, to prove $1$-regular convergence for
%quadrangulations, it suffices to prove that any simple loop of
%diameter $\leq \delta$ and made of edges of $S_\bq$, splits $S_\bq$
%into two connected components, one of which is of diameter $\leq
%\eps$.

\subsection{Tree encoding of quadrangulations}\label{sec:coding-with-trees}

We briefly introduce the second main ingredient, the Schaeffer
bijection. Let $\bT_n$ be the set of pairs $(\bt,\bl)$ where $\bt$ is
a rooted planar tree with $n$ edges, and $\bl$ is a function from the
set of vertices of $\bt$ to $\N=\{1,2,\ldots\}$, such that
$|\bl(x)-\bl(y)|\leq 1$ if $x$ and $y$ are neighbors. Then the set
$\bQ_n$ is in one-to-one correspondence with $\bT_n$. More precisely,
this correspondence is such that given a graph representation of
$\bq\in \bQ_n$ on a surface, the corresponding $(\bt,\bl)\in \bT_n$
can be realized as a graph whose vertices are
$V(\bt)=V(\bq)\setminus\{x_*\}$, where $x_*$ is the origin vertex of
the root edge, and $\bl$ is the restriction to $V(\bt)$ of the
function $\bl(x)=\dgr_\bq(x,x_*),x\in V(\bq)$. Moreover, the edges of
$\bt$ and $\bq$ only intersect at vertices. The root vertex
of $\bt$ is the tip of the root edge of $\bq$, so it lies at
$\dgr_\bq$-distance $1$ from $x_*$.

Let $x(0)$ be the root vertex of $\bt$, and given
$\{x(0),\ldots,x(i)\}$, and let $x(i+1)$ be the first
child\footnote{For the natural order inherited from the planar
  structure of $\bt$} of $x(i)$ not in $\{x(0),\ldots,x(i)\}$ if there
is any, or the parent of $x(i)$ if there is not. This procedure stops
at $i=2n$, where we are back to the root and have explored all
vertices of the tree.  We let $C_i=\dgr_\bt(x(i),x(0))$, and
$L_i=\bl(x(i))$. Both $C$ and $L$ are extended by linear interpolation
between integer times into continuous functions, still called $C,L$,
with duration $2n$. The contour process $C$ of $\bt$ is the usual
Harris walk encoding of the rooted tree $\bt$, and the pair $(C,L)$
determines $(\bt,\bl)$ completely. In the sequel, we will use the fact
that $x(i)$ can be identified with a vertex of $\bq$.

A simple consequence (see \cite[Lemma 3.1]{legall06}) of the
construction is that for $i<j$, 
\begin{equation}\label{eq:2}
\dgr_\bq(x(i),x(j))\leq L_i+L_j-2\min_{i\leq k\leq j}L_k+2\, .
\end{equation}

\subsection{Estimates on the lengths of 
geodesics}\label{sec:estim-lengths-geod}

Our last ingredient is a slight rewriting of the estimates of Le Gall
\cite{legall06} on geodesic paths in quadrangulations in terms of
encoding processes. Precisely, let $C^n,L^n$ be the contour and label
process of a uniform random element $\bt_n$ of $\mathbb{T}_n$, and let
$\bq_n$ be the quadrangulation that is the image of this element by
Schaeffer's bijection. In particular, $\bq_n$ is a random uniform
element of $\bQ_n$. Also, recall that a graphical representation
$\TT_n$ of $\bt_n$ can be drawn on the representation $\s_{\bq_n}$ of
Sect.\ \ref{sec:turn-quadr-into}, in such a way that the vertices of
$\TT_n$ are $\v_{\bq_n}\setminus\{x_*\}$, where $x_*$ is the root
vertex, and $\TT_n$ intersects edges $\e_{\bq_n}$ only at vertices.
For simplicity we let $V_n=V(\bq_n)$, $\dgr_n=\dgr_{\bq_n}$,
$\s_n=\s_{\bq_n}$, $d_n=d_{\bq_n}$.

The main result of \cite{legall05} says that the convergence in
distribution in ${\cal C}([0,1],\R)^2$ holds:
\begin{equation}\label{eq:1}\left(\left(\frac{1}{\sqrt{2n}}
C^n_{2nt}\right)_{0\leq t\leq 1},
\left(\left(\frac{9}{8n}\right)^{1/4}L^n_{2nt}\right)_{0\leq t\leq
  1}\right) \build\longrightarrow_{n\to\infty}^{(d)}(\ome,\oz)\, ,
\end{equation} where
$(\ome,\oz)$ is the Brownian snake conditioned to be positive
introduced by Le Gall \& Weill \cite{legweill}.  Moreover, it is shown
in \cite{legall06} that the laws of $[V_n,n^{-1/4}\dgr_n]$ form a
relatively compact family in the set of probability measures on $\M$
endowed with the weak topology. Since $V_n$ is $3$-dense in $\s_n$,
the same holds for $[\s_n,n^{-1/4}d_n]$. We argue as in
\cite{legall06}, and assume by Skorokhod's theorem that the trees
$\bt_n$ (hence also the quadrangulations $\bq_n$) are defined on a
common probability space on which we have, almost-surely
\begin{itemize}
\item
$[\s_n,n^{-1/4}d_n]\to [S,d]$, some random limiting space in $\pmc$,
  along some subsequence $n_k\to \infty$, and 
\item
the convergence (\ref{eq:1}) holds a.s.\ along this subsequence.
\end{itemize}
From this point on, we will always assume that $n$ is taken along this
subsequence. In particular, we have that $\diam S=\lim_n n^{-1/4}\diam
\s_n\geq \lim_n\sup n^{-1/4}L^n=\sup \ov{Z}>0$ a.s., so $S$ is not
reduced to a point and Theorem \ref{sec:grom-hausd-conv-1} may be
applied if we check that the convergence is $1$-regular.  We are going
to rely on proposition 4.2 of \cite{legall06}, which can be rephrased
as follows.

\begin{prp}\label{sec:estim-lengths-geod-2}
The following property is true with probability $1$. Let $i_n,j_n$ be
integers such that $i_n/2n\to s,j_n/2n\to t$ in $[0,1]$, where $s,t$
satisfy
$$\ome_s=\inf_{s\wedge t\leq u\leq s\vee
  t}\ome_u<\ome_t\, .$$ For $n\geq 1$, let $\gamma_n$ be
a path in $\bq_n$ between $x(i_n)$ and $x(j_n)$ with the notation of
Sect.\ \ref{sec:coding-with-trees}. Then it holds that
$$\liminf_{n\to\infty}n^{-1/4}{\rm length}(\gamma_n)>0\, .$$
\end{prp}

In \cite{legall06}, this proposition was a first step in the proof of
the fact that a limit in distribution of $(V_n,\dgr_n)$ can be
expressed as a quotient of the continuum tree with contour function
$\ome$: this lemma says that two points of the latter such that one is
an ancestor of the other are not identified. Le Gall completed this
study by exactly characterizing which are the points that are
identified. 

\section{Proof of Theorem \ref{T1}}\label{sec:proof-theorem-reft1-1}

\begin{lmm}\label{sec:proof-theorem-reft1}
Almost-surely, for every $\eps>0$, there exists a $\delta\in(0,\eps)$
such that for $n$ large enough, any simple loop $\gamma_n$ made of
edges of $\s_n$, with diameter $\leq n^{1/4}\delta$, splits $\s_n$ in
two Jordan domains, one of which has diameter $\leq n^{1/4}\eps$.
\end{lmm}

\proof We argue by contradiction, assuming there exist simple loops
$\gamma_n$ made of edges of $\s_n$, with diameters $o(n^{1/4})$ as
$n\to\infty$, such that the two Jordan domains bounded by $\gamma_n$
are of diameters $\geq n^{1/4}\eps$, where $\eps>0$ is some fixed
constant.  Let $l_n$ be the minimal label on $\gamma_n$, i.e.\ the
distance from the root vertex $x_*$ to $\gamma_n$. Then all the labels
of vertices that are in a connected component $D_n$ of $\s_n\setminus
\gamma_n$ not containing $x_*$ are all larger than $l_n$, since a
geodesic from $x_*$ to any such vertex must pass through
$\gamma_n$. 

\begin{figure}[!h]
\begin{center}
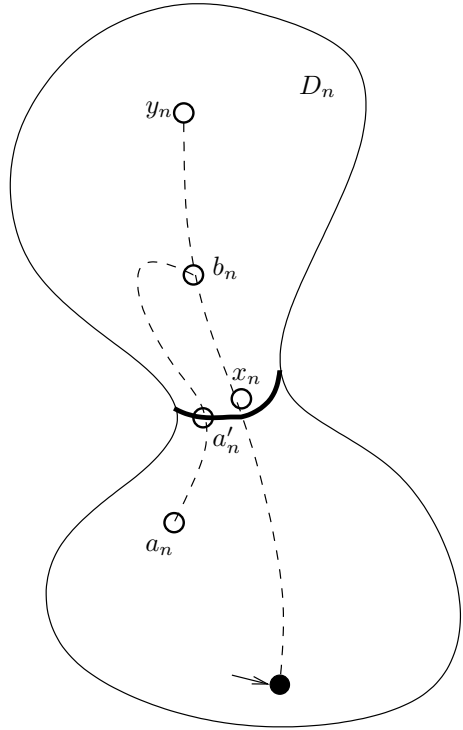
\caption{Illustration of the proof. The surface $\s_n$ is depicted as
  a sphere with a bottleneck circled by $\gamma_n$ (thick line). The
  root edge of the quadrangulation is drawn at the bottom, and the
  tree $\TT_n$ originates from its tip. In dashed lines are
  represented the two branches of $\TT_n$ that are useful in the
  proof: one enters the component $D_n$, and the other goes out after
  entering, identifying strict ancestors in the limit} \label{fig:1}
\end{center}
\end{figure}

The intuitive idea of the proof is the following. Starting from the
root of the tree $\TT_n$, follow a maximal simple path in $\TT_n$ that
enters in $D_n$ at some stage. If all such paths remained in $D_n$
after entering, then all the descendents of the entrance vertices
would have labels larger than that of the entrance vertex, a property
of zero probability under the limiting Brownian snake measure, see
\cite[Lemma 2.2]{legall06} and Lemma \ref{sec:estim-lengths-geod-1}
below. Thus, some of these paths must go out of $D_n$ after entering,
but they can do it only by passing through $\gamma_n$, which entails
that strict ancestors in $\TT_n$ will be at distance $o(n^{1/4})$, and
this is prohibited by Proposition \ref{sec:estim-lengths-geod-2}. This
is summed up in Figure \ref{fig:1}, which gathers some of the
notations to come.

We proceed to the rigorous proof. Take a vertex $y_n$ in $D_n$. As a
vertex of $\TT_n$, it can be written in the form $x(j_n)$ for some
$j_n$. Let $j'_n$ be the first integer $j\geq j_n$ such that $x(j)$ is
at $d_n$-distance $\leq 1$ from $\gamma_n$. Such a $j$ exists because
of the way edges of $\TT_n$ are drawn (entailing that the ancestral
path in $\TT_n$ from $x_n$ to the root of $\TT_n$ must itself pass at
distance $\leq 1$ from $\gamma_n$, since the root of $\TT_n$ is at
distance $1$ from $x_*$ and $x_n$ lies in $D_n$) and the label
$\bl(x(j'_n))$ is at most $\max_{z\in \gamma_n}
\bl(z)+1=l_n+o(n^{1/4})$. Moreover, for $k\in[j_n,j'_n]$, the vertex
$x(k)$ is in $D_n$, and in particular, its label is $\geq
l_n$. Applying the bound (\ref{eq:2}) to the times $j_n,j'_n$, we get
that $d_n(y_n,x(j'_n))\leq \bl(y_n)+l_n -2l_n +o(n^{1/4})$.  Since by
hypothesis the diameter of $D_n$ is at least $n^{1/4}\eps$, it is thus
possible to choose $y_n$ with label $\bl(y_n)\geq l_n+n^{1/4}\eps/2$.

We let $x_n$ be the first ancestor of $y_n$ in $\TT_n$ lying at
$d_n$-distance $\leq 1$ from $\gamma_n$, so that
$\bl(x_n)=l_n+o(n^{1/4})$.  Take $i_n<j_n$ such that $i_n$ is a time
encoding $x_n$, so that $C^n_{i_n}=\inf_{i_n\leq r\leq j_n} C^n_r$. Up
to further extraction, we may and will assume that
$$(9/8n)^{1/4}l_n\to l\, , \qquad i_n/2n\to s\, ,\qquad j_n/2n\to t\,
.$$ Then $s\leq t$ and $\ome_s\leq \ome_u$ for $u\in[s,t]$. More
precisely, we have $\oz_s=l$ and $\oz_t\geq l+(9/8)^{1/4}\eps/2$,
which implies $s<t$, and $\ome_s<\ome_t$. In terms of the continuum
tree encoded by $\ome$, this amounts to the fact that $s,t$ encode two
vertices such that the first is an ancestor of the second, and that
are not the same because the snake $\oz$ takes distinct values at
these points. We will need the following technical statement:

\begin{lmm}\label{sec:estim-lengths-geod-1}
Assume that $s>0$.  With probability $1$, there exist $\eta>0$ and
integers $i'_n,k_n,r_n$ with $i_n\leq i_n'<k_n<r_n<j_n$ such that
$i_n'/2n\to s'\in [s,t)$, that satisfy for $n$ large enough:
$$C^n_{i'_n}=\inf_{i'_n\leq r\leq j_n}C^n_r\, ,\qquad
C^n_{r_n}=\inf_{k_n\leq r\leq j_n}C^n_r\, ,$$ 
and 
$$(2n)^{-1/2}C^n_{i'_n}\to \ome_{s'}=\ome_s\, ,$$ so that
$i_n',k_n,r_n$ encode vertices $x_n',a_n,b_n$ of $\TT_n$ such that
$x_n\prec x_n'\prec b_n\prec y_n$ and $b_n\prec a_n$, where $\prec$
denotes ``is an ancestor of''. Moreover, $d_n(x_n,x'_n)=o(n^{1/4})$
and the labels satisfy
$$\bl(b_n)=L^n_{r_n}>L^n_{i_n}+\eta n^{1/4}=\bl(x_n)+\eta n^{1/4}\, ,$$
and
$$\bl(a_n)=L^n_{k_n}\leq L^n_{i_n}-\eta n^{1/4}=\bl(x_n)-\eta
n^{1/4}\, .$$
\end{lmm} 

The statement says roughly the following: there exist subtrees of
$\TT_n$ branching on a vertex $b_n$ of the ancestral line from $x_n$
to $y_n$ that attain labels that are significantly smaller (in the
scale $n^{1/4}$) than $\bl(y_n)$, but such that $\bl(b_n)$ is
significantly larger than $\bl(y_n)$.  As mentioned in the
Introduction, a rigorous proof needs some prerequisites on continuum
trees, and is really a re-writing of the proof of \cite[Proposition
  4.2, pp.649--650]{legall06}, using the fact \cite[Lemma
  2.2]{legall06} that the positive Brownian snake has no increase
points except $0$. We only explain in detail the role of $i'_n$ in the
statement. Introducing $i'_n\neq i_n$ may be necessary if it happens
that two macroscopic subtrees branch just above $x_n$. This happens if
$\ome$ attains a local minimum equal to $\ome_s$ at a time
$s'\in(s,t)$, which is the most pathological situation that can occur
since local minima of a Brownian motion are pairwise distinct and
realized only once. In this case, we take $i'_n$ so that $i'_n/2n\to
s'$ and $C^n_{2n\cdot}$ achieves a local minimum at $i'_n$. Then, the
vertex encoded by $i'_n$ is encoded by another time $i''_n<i'_n$ such
that $i''_n/2n\to s$, which together with (\ref{eq:2}) implies the
property $d_n(x_n,x'_n)=o(n^{1/4})$. If $\ome_s<\ome_u$ for every
$u\in (s,t]$, we simply take $i'_n=i_n$.

\medskip

Now back to the proof of Lemma \ref{sec:proof-theorem-reft1}. Because
of the property of the label of $a_n$, it does not lie in $D_n$,
however, its ancestor $b_n$ does because it is on the ancestral path
from $x_n$ to $y_n$. Hence some ancestor of $a_n$ must belong to
$\gamma_n$, and let $a'_n$ be the youngest of these (the highest in
the tree), and take $k'_n\in(k_n,r_n)$ encoding $a_n'$. Since $x_n$ is
at distance at most $1$ from $\gamma_n$, we obtain that
$d_n(a'_n,x_n)=o(n^{1/4})$. However, if $k'_n/2n\to v,r_n/2n\to u$,
taking again an extraction if necessary, then we have $\ome_s<
\ome_u\leq\ome_v$ because of the ancestral relations $C^n_{i_n}\leq
C^n_{r_n}\leq C^n_{k'_n}$, and the fact
$$\oz_s=\lim_{n\to\infty} (9/8n)^{1/4}L^n_{i_n}= l<l+(9/8)^{1/4}\eta\leq
\lim_{n\to\infty} (9/8n)^{1/4}L^n_{r_n}=\oz_u\, .$$ Now the statements
$d_n(a'_n,x_n)=o(n^{1/4})$ and $\ome_s<\ome_v$ together can only hold
with zero probability by Proposition \ref{sec:estim-lengths-geod-2}.

It remains to rule out the possibility that $s=0$, i.e.\ that
$\gamma_n$ lies at distance $o(n^{1/4})$ from $x_*$. To see that this
is not possible, argue as in the beginning of the proof and take
$x_n,y_n$ respectively in the two disjoint connected components of
$\s_n\setminus \gamma_n$, and with labels $\bl(x_n)\wedge \bl(y_n)\geq
n^{1/4}\eps/2$. By symmetry, assume that $x_n=x(i_n)$ and $y_n=x(j_n)$
with $i_n<j_n$. Now take the least integer $k_n\in[i_n,j_n]$ such that
$x(k)$ belongs to $\gamma_n$. Such a $k$ has to exist because any path
from $x_n$ to $y_n$ in $\s_n$ must pass through $\gamma_n$. Then
$L^n_{k_n}=\bl(x(k_n))=o(n^{1/4})$. Up to extraction, assume
$i_n/2n\to s,k_n/2n\to u,j_n/2n\to t$. Then $\oz_u=0<\oz_s\wedge
\oz_t$, so that $s<u<t$, and this contradicts the fact that $\oz$ is
strictly positive on $(0,1)$, which is a consequence of
\cite[Proposition 2.5]{legweill}. \cq

\medskip

We claim that this lemma is enough to obtain $1$-regularity of the
convergence, and hence to conclude by Theorem
\ref{sec:grom-hausd-conv-1} that the limit $(S,d)$ is a sphere. First
choose $\eps<\diam S/3$ to avoid trivialities. Let $\gamma_n$ be a
loop in $\s_n$ with diameter $\leq n^{1/4}\delta$. The boundary of the
union of the closures of faces of $\bq_n$ that are hit by $\gamma_n$
is made of pairwise disjoint simple loops of edges of $\s_n$. If $x,y$
are elements of this union of faces, and since a face of $\s_n$ has
diameter less than $3$, there exist points $x',y'$ of $\gamma_n$ at
distance at most $3$ from $x,y$ respectively, so that the diameters of
these loops all are $\leq n^{1/4}\delta+6$. By the Jordan Curve
Theorem, each of these loops splits $\s_n$ into two simply connected
components, one of which has diameter $\leq n^{1/4}\eps$, and one of
which contains $\gamma_n$ entirely. It suffices to justify that these
two properties (being of diameter $\leq n^{1/4}\eps$ and containing
$\gamma_n$) hold simultaneously for some loop in the family to
conclude that $\gamma_n$ is homotopic to $0$ in its
$\eps$-neighborhood. So assume the contrary: the component not
containing $\gamma_n$ of every loop is of diameter $\leq
n^{1/4}\eps$. By definition, any point in the complement of the union
of these components is at distance at most $3$ from some point of
$\gamma_n$. Take $x,y$ such that $d_n(x,y)=\diam(\s_n)$. Then there
exist points $x',y'$ in $\gamma_n$ at distance at most $n^{1/4}\eps+3$
respectively from $x,y$, and we conclude that $d_n(x',y')\geq
\diam(\s_n)-6-2n^{1/4}\eps>n^{1/4}\delta$ for $n$ large enough by our
choice of $\eps$, a contradiction.

\def\polhk#1{\setbox0=\hbox{#1}{\ooalign{\hidewidth
  \lower1.5ex\hbox{`}\hidewidth\crcr\unhbox0}}}

\end{document}